\documentclass{SCAE}
\numberwithin{equation}{section}
\begin{document}
\newcommand\med{\medskip}
\newcommand\eq{\eqno}
\def\m{\mu}
\def \Cal{\cal}
\newcommand{\IR}{\hbox{\rm l\negthinspace R}}
\def\tag {\eqno}
\newcommand {\ty}{\overline y}
\newcommand {\tz}{\overline z}

\newcommand{\w}{\omega}
\newcommand{\CW}{\Omega}
\newcommand{\E}{\mathbb{E}}
\newcommand{\CF}{{\cal F}}
\newcommand{\ld}{\lambda}
\newcommand{\CR}{\mathbb{R}}
\newcommand{\CE}{\mathbb{E}}
\newcommand{\CV}{\mathbb{V}}
\newcommand{\Ce}{{\cal E}}
\newcommand{\CP}{{\cal P}}

\Year{2016} %
\Month{May}
\Vol{59} %
\No{5} %
\BeginPage{945} %
\EndPage{954} %
\AuthorMark{CHEN Z J}
\ReceivedDay{February 21, 2015}
\AcceptedDay{October 12, 2015}
\PublishedOnlineDay{; published online December 2, 2015}
\DOI{10.1007/s11425-015-5095-0} 

\title{Strong laws of large numbers for sub-linear expectations}{}


\author[]{CHEN ZengJing}{}

\address[]{School of Mathematics, Shandong University, Jinan {\rm 250100}, China;}

\Emails{zjchen@sdu.edu.cn;}\maketitle


 {\begin{center}
\parbox{14.5cm}{\begin{abstract}
We investigate three kinds of strong laws of large numbers for capacities with a new notion of
independently and identically distributed (IID) random variables for sub-linear expectations initiated by Peng.
It turns out that these theorems are natural and fairly neat extensions of the classical Kolmogorov¡¯s strong law
of large numbers to the case where probability measures are no longer additive. An important feature of these
strong laws of large numbers is to provide a frequentist perspective on capacities.\vspace{-3mm}
\end{abstract}}\end{center}}

 \keywords{capacity, strong law of large numbers, independently and identically distributed, nonlinear expectation}

 \MSC{60F15, 60G50}

\renewcommand{\baselinestretch}{1.2}
\begin{center} \renewcommand{\arraystretch}{1.5}
{\begin{tabular}{lp{0.8\textwidth}} \hline \scriptsize
{\bf Citation:}\!\!\!\!&\scriptsize     Chen Z J. Strong laws of large numbers for sub-linear expectations. Sci China Math, 2016, 59: 945¨C954,
doi: 10.1007/s11425-015-5095-0\vspace{1mm}
\\
\hline
\end{tabular}}\end{center}

\baselineskip 11pt\parindent=10.8pt  \wuhao
\section{Introduction}
The classical strong laws of large numbers (strong LLN) as fundamental limit theorems in probability
theory play a fruitful role in the development of probability theory and its applications. The key in
the proofs of these limit theorems is the additivity of probability measures and mathematical expectations.
However, such an additivity assumption is not feasible in many areas of applications because
many uncertain phenomena cannot be well modelled using additive probabilities or additive expectations.
More specifically, motivated by some problems in mathematical economics, statistics, quantum mechanics
and finance, a number of papers have used non-additive probabilities (called capacities) and nonlinear
expectations (for example Choquet integral/expectation, $g$-expectation) to describe and interpret the
phenomena which are generally nonadditive (see \cite{chenepstein,Feynman,Gilboa,Huber,Peng:1997,Peng:1999,Wakker,Wasserman}). A natural question is what is
the law of large numbers under nonadditive probabilities or nonlinear expectations? Recently, motivated
by the risk measures, super-hedge pricing and model uncertainty in finance, Peng\cite{peng1,peng2,peng3,peng4,peng5,peng6} initiated the
notion of independently and identically distributed (IID) random variables under sub-linear expectations.
Under this framework, he proved a weak law of large numbers (LLN) and a central limit theorem (CLT).
In this paper, we investigate three strong laws of large numbers for capacities in Peng¡¯s framework. All
of them are natural and fairly neat extensions of the classical Kolmogorov¡¯s strong law of large numbers,
but the proofs here are different from the original proofs of the classical strong law of large numbers.

Now we describe the problem in more details. For a given set $\CP$
of multiple prior probability measures on $(\Omega,\CF)$, we
define a pair $(\CV, v)$ of capacities by
$$\CV(A):=\sup_{P\in \CP}P(A),\;\quad  v(A):=\inf_{P\in\CP}P(A), \quad \forall A\in\CF.$$
The corresponding Choquet integrals/expecations $(C_{\CV},C_v)$
are defined by
$$
C_V[X]:=\int_0^\infty V(X\ge t)dt+\int_{-\infty}^0[V(X\ge t)-1]dt,
$$
where $V$ is replaced by $\CV$ and $v$, respectively.

The pair of so-called maximum-minimum expectations $(\CE,\Ce)$ are
defined by
$$\CE[\xi]:=\sup_{P\in\CP}E_P[\xi], \quad \Ce[\xi]:=\inf_{P\in\CP}
E_{P}[\xi].$$ Here and in the sequel, $E_P$ denotes  the classical
expectation under probability $P.$

In general, the relation between Choquet integral and maximum-minimum expectations are as follows: For any random variable $X$,
$$
\CE[X]\leq C_\CV[X],\quad\quad C_v[X]\leq \Ce[X]
$$
Note that under some very special assumptions on $\CP$ and $\CV$,
both inequalities could become equalities (see\cite{Gilboa,Huber,Schmeidler}).

Given a sequence $\{X_i\}_{i=1}^\infty$ of IID random variables for
capacities, the earlier papers related to strong laws of large
numbers for  capacities can be found in\cite{Dow,Walley}. However, the more general results for strong
laws of large numbers for capacities were given by Maccheroni and Marinacci \cite{Marinacci2}, Marinacci \cite { Marinacci1}
and Epstein and Schneider \cite{epstein}. They show that, on full set, any
cluster point of empirical averages lies between the lower Choquet
integral $C_v[X_1]$ and the upper Choquet integral $C_\CV[X_1]$ with
probability one under capacity $v$, i.e.,
$$v\left(\omega\in\Omega:
C_v[X_1]\leq\liminf\limits_{n\rightarrow\infty}\frac
1n\sum\limits_{i=1}^n X_i(\omega)
\leq\limsup\limits_{n\rightarrow\infty}\frac1n\sum\limits_{i=1}^nX_i(\omega)\leq
C_{\CV}[X_1]\right)=1.$$ Marinacci \cite{ Marinacci1} obtains his result under
the assumptions that $\CV$ is a totally monotone capacity on a
Polish space $\Omega$, random variables $\{X_i\}_{i=1}^{\infty}$ are
bounded or continuous. Epstein and Schneider\cite{epstein}  also shows the same
result under the assumptions that $\CV$ is rectangular and  the set
$\CP$ is finite.

Since the gap  between the Choquet integrals $C_{\mathbb{V}}[X]$  and $C_v[X]$ is bigger than that of the maximum-minimum expectations
$ \CE[X]$ and $\Ce[X]$ for all $X$, it is of interest to ask
whether we can obtain a more precise result if the Choquet
integrals/expctations in the above equality are replaced by
maximum-minimum expectations, i.e.,
$$
v\left(\omega\in \Omega:
\Ce[X_1]\leq\liminf\limits_{n\to\infty}\frac1n\sum\limits_{i=1}^n
X_i(\omega) \leq\limsup\limits_{n\to\infty}\frac
1n\sum\limits_{i=1}^nX_i(\omega)\leq \CE[X_1]\right)=1.
$$

The first result in this paper is to show that the above equality
is still true in Peng's framework. Furthermore, motivated by this result, we establish
two new laws of large numbers. The first is to show that there
exist two cluster points of empirical averages which reach the
minimum expectation $\Ce[X_1]$ and the maximum expectation
$\CE[X_1]$, respectively under capacity $\CV$, i.e.,
$$
\CV\left(\omega\in \Omega: \limsup
\limits_{n\to\infty}\frac1n\sum\limits_{i=1}^nX_i(\omega)=
\CE[X_1]\right)=1,$$
$$
\CV\left(\omega\in
\Omega:\liminf\limits_{n\to\infty}\frac1n\sum\limits_{i=1}^n
X_i(\omega)= \Ce[X_1]\right)=1.$$

The second  is to prove that the
cluster set of empirical averages is coincide with the interval
between minimum expectation $\Ce[X_1]$ and maximum expectation
$\CE[X_1],$ i.e., let $C(\{x_n\})$ be the cluster set of
$\{x_n\},$ then, for any $b\in \left[\Ce[X_1],
\CE[X_1]\,\right],$
$$
\CV\left( \omega\in\Omega: b\in  C\left(\left\{\frac1n\sum\limits_{i=1}^n
X_i(\omega)\right\}\right)\right)=1.
$$
Obviously, if either $\CV$ or $v$ in the above results is a probability measure, all of our main results are natural and
fairly neat extensions of the classical Kolmogorov's strong law of
large numbers. Moreover, an important feature of our strong laws
of large numbers is to provide a frequentist perspective on
capacities.

\section{Notation and lemmas}

In order to prove our results in Peng's framework, we shall recall briefly the notions of  both IID random
variables and sub-linear expectations initiated by Peng\cite{peng3}.

Let $(\CW,\CF)$ be a measurable space, and $\mathcal {L}$ be a subset
of  all random variables on $(\CW,\CF)$ such that for any $A\in\CF,$ $I_A\in \mathcal {L},$ where $I_A$ is the indictor function of event $A$.

\begin{definition}
A functional $\mathbb{E}$ on $\mathcal {L}
\mapsto(-\infty,+\infty)$ is called a sub-linear expectation, if
it satisfies the following properties: for all $X,Y\in\mathcal
{L},$

 (a) Monotonicity: $X\geq Y$ implies
$\mathbb{E}[X]\geq\mathbb{E}[Y]$.

 (b) Constant preserving:
$\mathbb{E}[c]=c,\forall c\in\mathbb{R}$.

 (c) Sub-additivity:
$\mathbb{E}[X+Y]\leq\mathbb{E}[X]+\mathbb{E}[Y]$.

 (d) Positive homogeneity: $\mathbb{E}[\lambda
X]=\lambda\mathbb{E}[X],\forall\lambda\geq0$.

\end{definition}

Given a sub-linear expectation $\CE$, let us denote the conjugate
expectation $\Ce$ of  sub-linear $\CE$  by
$$
\Ce[X]:=-\CE[-X], \quad  \forall X\in \mathcal {L}
$$
Obviously, for all $X\in\mathcal {L},$ $\Ce[X]\leq \CE[X].$
By  the sub-additivity of $\CE$, we have the following lemma.
\begin{lemma}\label{le1}
If $X, Y\in \mathcal {L},$ then
$$\Ce[X]\leq \CE[X+Y]-\CE[Y].
$$
\end{lemma}
Given a sub-linear expectation, we can define   a pair of capacities $(\CV, v)$  as follows:
\begin{definition}\label{definition0} A pair $(\CV,v)$ of capacities is said to be generated by a sub-linear expectation $\CE,$ if
$$
\CV(A):=\CE[I_A], \quad  v(A):= \Ce[I_A],\quad \forall A\in\CF.
$$
\end{definition}
It is easy to check that such capacities have the following properties:
\begin{lemma}\label{definition1}
(1) $\CV(\varnothing)=v(\varnothing)=0, \CV(\Omega)=v(\Omega)=1$.

(2) $\CV(A)\leq \CV(B),$ $v(A)\le v(B)$,
whenever $A\subset B$ and $A,B\in\CF$.

(3)
$\CV(A\cup B)\leq  \CV(A) +\CV(B)$, $A,B\in\CF$.

(4)
$
\CV(A)+v(A^c)=1,\quad \forall A\in\CF,
$
where $A^c$ is the complement set of $A.$

\end{lemma}

Motivated by the notion of IID random variables under sub-linear expectations initiated by Peng \cite{peng3}, we adopt the following  notion of IID random variables under sub-linear expectations to study strong law of large numbers for non-additive probabilities.

\begin{definition}\label{de3}

{\bf Independence.} Suppose that $Y_1,Y_2,\cdots,Y_n$ is a
sequence of random variables such that $Y_i \in\mathcal {L}.$
Random variable $Y_n$ is said to be independent to
$X:=(Y_1,\cdots,Y_{n-1})$ under $\mathbb{E}$, if for each Borel-measurable
 function $\varphi$ on $\mathbb{R}^n$ with $\varphi(X,Y_n)\in
\mathcal {L}$ and $\varphi(x,Y_n)\in\mathcal {L}$ for each $x\in
\mathbb{R}^{n-1},$ we have
$$\mathbb{E}[\varphi(X,Y_n)]=\mathbb{E}[\overline{\varphi}(X)],$$
where $\overline{\varphi}(x):=\mathbb{E}[\varphi(x,Y_n)]$ and
$\overline{\varphi}(X)\in\mathcal {L}$.

\noindent{\bf Identical distribution.} Random variables $X$ and $Y$
are said to be identically distributed, denoted by
$X\overset{d}{=}Y$, if for each Borel-measurable
 function $\varphi$ such that $\varphi(X),
\; \varphi(Y)\in \mathcal {L}$,
$$\mathbb{E}[\varphi(X)]=\mathbb{E}[\varphi(Y)].$$

\noindent{\bf IID random variables.} A sequence of random variables
$\{X_i\}_{i=1}^\infty$ is said to be IID, if
$X_i\overset{d}{=}X_1$ and $X_{i+1}$ is independent to
$Y:=(X_1,\cdots,X_i)$ for each $i\ge 1.$

\end{definition}

The following lemma shows the relation between our  independence
and pairwise independence in \cite{Marinacci1}.

\begin{lemma}\label{le2}
Suppose that $X,Y \in \mathcal {L}$ are two random variables.
$\CE$ is a sub-linear expectation and $(\CV,v)$ is the pair of
capacities generated by  $\CE.$ If  random variable $X$ is
independent to $Y$ under $\CE$, then $X$ is also
independent to $Y$ under capacities $\CV$ and $v,$ i.e. for all
subsets $D$ and $G\in \mathcal{B}(\CR),$
$$
V(X\in D, Y\in G)=V(X\in D)V(Y\in G)
$$
holds for both capacity $\CV$ and $v$.
\end{lemma}

\begin{proof} If we choose $\varphi(x,y)=xy,$
by the definition of  independence in Definition \ref{de3}, it is easy to obtain the independence for events,
$$
\mathbb{V}(X\in D, Y\in G)=\CE[I_{\{X\in D\}}I_{\{Y\in G\}}]=\CE[\varphi(I_{\{X\in D\}}, I_{\{Y\in G\}})]=\mathbb{V}(X\in D)\mathbb{V}(Y\in G).
$$
Similarly, we can prove that $X$ is independent to $Y$ under capacities $v$
by choosing $\varphi(x,y)=-xy$.
\end{proof}

Chen et al. \cite{chenwu} prove that Borel-Cantelli Lemma is still true for capacity under some
assumptions.

\begin{lemma}\label{le3}(See [2, Lemma 2.2]).
Let $\{A_n,n\geq1\}$ be a sequence of events in $\CF$ and $(\CV,v)$
be a pair of capacities generated by sub-linear expectation $\CE$.

(1) If $\sum\limits_{n=1}^\infty\CV(A_n)<\infty,$ then $\CV\left(\bigcap\limits_{n=1}^\infty\bigcup\limits_{i=n}^\infty
 A_i\right)=0.$

(2) If further $\mathbb{V}$ is upper continuous and $\{A_n^c\}_{n=1}^{\infty}$ are mutually independent with respect to $v$, i.e., for any $n\in\mathbb{N}$, $$v\left(\bigcap\limits_{i=n}^\infty
A_i^c\right)=\prod_{i=n}^\infty v(A_i^c).$$ If
$\sum\limits_{n=1}^\infty{\CV}(A_n)=\infty $, then
$$\CV\left(\bigcap\limits_{n=1}^\infty\bigcup\limits_{i=n}^\infty A_i\right)=1.$$

\end{lemma}

Suppose that $C_b(\CR)$  is the  set of  all  continuous and bounded functions on $ \CR$ and $C_b^2(\CR)$ is the set of all continuous and bounded  functions on $ \CR$  whose second derivatives exist in $C_b(\CR)$.

With the notion of IID under sub-linear expectation,  we can obtain the following Lemma.

\begin{lemma}\label{le4}
Let $\{X_i\}_{i=1}^\infty$ be a sequence of IID random
variables with finite means $\overline\mu:=\CE[X_1], \;\;
\underline \mu:=\Ce[X_1],$  and  $S_n:=\sum_{i=1}^n X_i$ with $S_0:=0.$
 Suppose $\CE[|X_1|^{1+\alpha}]
<\infty$ for some $\alpha>0.$  Then for any positive
function $\varphi \in C_b(\mathbb{R}),$
$$
\liminf \limits_{n\to \infty}\;\CE\left[\varphi \left(\frac {S_n}n \right)\right]\ge
\sup_{\underline \mu\leq x\leq\overline \mu}\varphi(x).
$$
\end{lemma}

\begin{proof}  We turn the proof into three steps.
  Let $x^*$ is the maximal point of $\varphi$ over $[\underline \mu, \overline \mu].$

\noindent{\bf Step 1.} We first prove that  if $\{X_i\}_{i=1}^{\infty}$ is a IID sequence, then
$$
\CE\left[\varphi\left( \frac 1n \sum_{i=1}^n X_i\right)\right]-\varphi\left(x^*\right)\ge  n \inf\limits_{x\in\CR}\left\{\CE\left[\varphi\left(x+\frac{X_{n-m}-x^*}n\right)\right]
-\varphi\left(x\right)\right\}
$$
In fact,  set $T_k:= \frac 1n\sum_{i=1}^k X_i$ with $T_0=0, k=1,2, \cdots n,$ and $y:=\frac {x^*}n.$
\begin{eqnarray}\label{eq1-10}
\nonumber
\CE[\varphi\left(T_n \right)]-\varphi\left(x^*\right)&&=\CE[\varphi\left(T_n \right)]-\CE[\varphi( T_{n-1}+ y)]\\\nonumber
&& \hspace*{1em}+\CE[\varphi\left( T_{n-1}+ y\right)]
-\CE\left[\varphi\left( T_{n-2}+2y\right)\right]+\cdots\\\nonumber
&& \hspace*{1em}+\CE\left[\varphi\left(T_{n-m}+ m y\right)\right]-\CE\left[\varphi\left( T_{n-(m+1)}+ (m+1)y\right)\right]+\cdots\\\nonumber
&& \hspace*{1em}+\CE\left[\varphi\left( T_1+ (n-1)y\right)\right]-\CE\left[\varphi\left(n y\right)\right]\\
&&=\sum\limits_{m=0}^{n-1}\left\{
\CE\left[\varphi\left( T_{n-m}+ m y\right)\right]-\CE\left[\varphi\left( T_{n-(m+1)}+(m+1) y\right)\right]\right\}.
\end{eqnarray}
We now evaluate each term inside the summation.
Let
$$h(x):=\CE\left[\varphi\left(x+\frac {X_{n-m}}n\right)\right].$$
  Then  because of independence  of $\{X_i\}_{i=1}^n,$
\begin{eqnarray}\label{eq8}
\nonumber \CE\left[\varphi \left( T_{n-m}+ m y \right)\right]
&=& \CE\left[\CE\left[\varphi\left( x+\frac {X_{n-m}}n \right)\right]\Big|_{x=T_{n-(m+1)}+ m y}\right]\\\nonumber
&=&\CE\left[h \left( T_{n-(m+1)}+ m y \right)\right].
\end{eqnarray}
Then by the sub-linearity of $\CE$ in  Lemma \ref{le1}, we have
\begin{eqnarray}\label{eq0-01}
\nonumber
&&\CE\left[\varphi\left( T_{n-m}+ m y\right)\right]-\CE\left[\varphi \left( T_{n-(m+1)}+(m+1) y\right)\right]\\\nonumber
&=&\CE\left[h\left( T_{n-(m+1)}+ my\right)\right]
-\CE\left[\varphi \left( T_{n-(m+1)}+m y+ y\right)\right]\\\nonumber
&\ge &\Ce\left[h\left( T_{n-(m+1)} +m y\right)
-\varphi \left( T_{n-(m+1)}+m y+y\right)\right]\\\nonumber
&\ge &\inf\limits_{x\in\CR}\left (h(x)-\varphi\left(x+y\right)\right)\\\nonumber
&=&\inf\limits_{x\in\CR}\left\{\CE\left[\varphi\left(x+\frac {X_{n-m}}n\right)\right]
-\varphi\left(x+\frac {x^*}n\right)\right\}\\\nonumber
&=&\inf\limits_{x\in\CR}\left\{\CE\left[\varphi\left(x+\frac {X_{n-m}-x^*}n\right)\right]
-\varphi\left(x\right)\right\}.
\end{eqnarray}
 It then follows that $\{X_i\}_{i=1}^{\infty}$ is identical. The proof of Step 1 is complete.

\noindent{\bf Step 2.}  For $\varphi\in C_b^2(\CR),$ We shall prove that
$$ \liminf\limits_{ n\to\infty}\; n\; \inf\limits_{x\in\CR}\left\{\CE\left[\varphi\left(x+\frac {X_{n-m}-x^*}n\right)\right]
-\varphi\left(x\right)\right\}\ge 0.$$
The Taylor expansion of function $\varphi$ implies that for some random variables $\{\theta_i\}_{i=1}^n$ valued in $[0, 1],$
\begin{eqnarray}\label{INE}
\varphi\left(x+\frac{X_i-x^*}n\right)-\varphi(x)
&=& \varphi'(x)\frac{X_i-x^*}n+J_n(x, X_i,x^*),
\end{eqnarray}
where $$J_n(x, X_i,x^*):=\left[\varphi'\left(x+\theta_i\frac{X_i-x^*}n\right)-\varphi'(x)\right]\frac{X_i-x^*}n,\quad 1\leq i\leq n.$$
Taking sub-linear expectation $\CE$ on both sides of  (\ref{INE}), and applying
the sub-linearity of $\CE$, we have
\begin{eqnarray}\label{EQ-1}\nonumber
-\CE[|J_n(x, X_i,x^*)|]+\CE\left[\varphi'(x)\frac{X_i-x^*}n\right]&\leq&
\CE\left[\varphi\left(x+\frac{X_i-x^*}n\right)-\varphi(x)\right].
\end{eqnarray}

Since $\CE[X_i]=\overline \mu,\;  $  $\Ce[X_i]=\underline\mu$ and $x^*\in [\underline \mu, \; \overline \mu]$,
$$\CE\left[\varphi'(x)\frac{X_i-x^*}n
\right]=(\varphi'(x))^+\; \frac {\overline\mu-x^*}n+(\varphi'(x))^-\;\frac {x^*-\underline\mu}n\ge 0.
$$
Therefore, we only need to prove that
\begin{equation}\label{1-1}
 \sum_{i=1}^n\sup_{x\in \CR}\CE[|J_n(x,X_i,x^*)|]\to 0, \quad n\to\infty.
 \end{equation}
In fact, for any $\epsilon >0,$ using H\"{o}lder and Chebyshev's inequalities and the fact that $\{X_i\}$ is identical, we get
\begin{eqnarray}\label{eq4h}\nonumber
&&\sum_{i=1}^n\sup_{x\in\CR}\CE[|J_n(x,X_i,x^*)|]\\\nonumber
&\leq& \sum_{i=1}^n \left\{\sup_{x\in\CR}\CE\left[
|J_n(x,X_i,x^*)| I_{\{|\frac{X_i-x^*}n |>  \epsilon\}}\right]+\sup_{x\in\CR}\CE\left[|J_n(x,X_i,x^*) |I_{\{|\frac{X_i-x^*}n|
\leq \varepsilon\}}\right]\right\}\\ \nonumber
&\leq& \sum_{i=1}^n\left\{\CE\left[\left(\sup_{x\in\CR}\left|\varphi'\left(x+\theta_i\frac{X_i-x^*}n\right)\right|+\sup_{x\in\CR}|\varphi'(x)|\right)\frac{|X_i-x^*|}n I_{\{|X_i-x^*|> n \epsilon\}}\right]\right.\\\nonumber
&&+\left.\CE\left[\sup_{x\in\CR}\left|\varphi''\left(x+\theta_i\bar{\theta}_i\frac{X_i-x^*}n\right)\right|\frac{(X_i-x^*)^2}{n^2}I_{\{|{X_i-x^*} |\leq n \epsilon\}}\right]\right\}\\\nonumber
&\leq& n\left\{\frac{2\|\varphi'\|}{n}(\mathbb{E}[|X_1-x^*|^{1+\alpha}])^{\frac{1}{1+\alpha}} (\mathbb{E}[I_{\{|X_1-x^*|> n \epsilon\}}])^{\frac{\alpha}{1+\alpha}}         +\frac{\epsilon}{n}\|\varphi''\|\mathbb{E}[|X_1-x^*|]\right\}\\\nonumber
&\leq& n\left\{\frac{2\|\varphi'\|}{n}(\mathbb{E}[|X_1-x^*|^{1+\alpha}])^{\frac{1}{1+\alpha}} \left(\frac{\mathbb{E}[|X_1-x^*|^{1+\alpha}]}{(n\epsilon)^{1+\alpha}}\right)^{\frac{\alpha}{1+\alpha}}  +\frac{\epsilon}{n}\|\varphi''\|\mathbb{E}[|X_1-x^*|]\right\}\\\nonumber
&\leq& n\left\{ \frac 2{n^{1+\alpha}\epsilon^\alpha}\|\varphi'\|\CE\left[|{X_1-x^*}|^{1+\alpha} \right] + \frac {\epsilon}{n}\|\varphi''\|\CE[|X_1-x^*| ]\right\}\\\nonumber
&=& \frac 2{(n\epsilon)^{\alpha}}\|\varphi'\|\CE\left[|{X_1-x^*}|^{1+\alpha} \right] + \epsilon\|\varphi''\|\CE[|X_1-x^*| ]\\\nonumber
&\to& \epsilon\|\varphi''\| \CE[|X_1-x^*|], \; \hbox{as }n\to \infty,
\end{eqnarray}
where $\{\bar{\theta}_i\}_{i=1}^{\infty}$ are random variables valued in $[0,1]$. For arbitrariness of $ \epsilon,$  we obtain the conclusion (\ref{1-1}).

Hence, the Lemma\ref{le4} hold for $\varphi\in C_b^2(\CR)$.

\noindent{\bf Step 3.}
If $\varphi\in C_b(\CR),$  then for any $\epsilon>0$  there exists $\overline \varphi\in C_b^2(\CR)$ such that $$
\sup_{x\in \CR}| \varphi(x)-\overline \varphi(x)|\leq \epsilon.$$

Apply Step 2 for function $\overline \varphi(x)$ and the fact that
\begin{eqnarray}\nonumber
&&\liminf_{n\to\infty}\CE\left[\varphi\left(\frac {S_n}n\right)\right]-\sup_{\underline  \mu\leq  x\leq \overline \mu}\varphi(x)\\\nonumber
&=&\liminf_{n\to\infty}\CE\left[\varphi\left(\frac {S_n}n\right)-\overline\varphi\left(\frac {S_n}n\right)+\overline\varphi\left(\frac {S_n}n\right)\right]-\sup_{\underline  \mu\leq  x\leq \overline \mu}[\varphi(x)-\overline\varphi(x)+\overline\varphi(x)]\\\nonumber
&\ge &\liminf_{n\to\infty}\CE\left[\overline\varphi\left(\frac {S_n}n\right)\right]-\sup_{\underline  \mu\leq  x\leq \overline \mu}\overline\varphi(x)-2\epsilon\\\nonumber
&\ge &-2\epsilon.
\end{eqnarray}
For arbitrariness of $ \epsilon,$  the proof of this Lemma is complete.
\end{proof}

\section{Main result }

The following theorem is our main result.
\begin{theorem}\label{th1}
Let $\{X_i\}_{i=1}^\infty$ be a sequence of IID random
variables for sublinear expectation $\CE$. Suppose
$\CE[|X_1|^{1+\alpha}]<\infty$ for some $\alpha\in (0, 1].$ Set
$\overline \mu:=\CE[X_1],$ $\underline \mu=\Ce[X_1]$  and $
S_n:=\sum\limits_{i=1}^n X_i.$ Then

(I)
\begin{equation}\label{M}
\CV\left(  \{\liminf\limits_{n\to\infty} S_n/n <\underline
\mu\}\bigcup \{\limsup\limits_{n\to\infty} S_n/n
>\overline \mu\}\right)=0.
\end{equation}
Also
\begin{equation}\label{M2}
v\left( \underline \mu\leq \liminf\limits_{n\to\infty} S_n/n \leq
\limsup\limits_{n\to\infty} S_n/n \leq \overline \mu\right)=1.
\end{equation}
If furthermore $\mathbb{V}$ is upper continuous, then

(II)\label{M1}
$$\CV\left( \limsup\limits_{n\to\infty}S_n/n= \overline
\mu\right)=1,\quad \CV\left( \liminf\limits_{n\to\infty} S_n/n=
\underline \mu\right)=1.
$$

(III) Suppose that  $C(\{x_n\})$ is the cluster set of a sequence of
 $\{x_n\}$ in $\CR,$ then, for any $b\in [\underline\mu, \overline
 \mu]$
 $$
 \CV\left (b\in C(\{S_n/n\})\right)=1.
 $$

\end{theorem}

\begin{proof}
(I) can be deduced  from [\cite{chenwu},Theorem 3.1] directly, we omit the details.

We now prove (II). If $\overline{\mu}=\underline{\mu}$, it is trivial.
Suppose $\overline{\mu}>\underline{\mu}$, then we only need to prove
that there exists an increasing subsequence $\{n_k\}$ of $\mathbb{N}$
such that for any $0<\epsilon< \overline \mu-\underline \mu,$
\begin{equation}\label{R2}
\CV\left(\bigcap_{m=1}^{\infty}\bigcup_{k=m}^{\infty}\{S_{n_k}/{n_k}\ge
\overline \mu-\epsilon\}\right)=1.
\end{equation}
Since $\mathbb{E}$ is upper continuous, we have
 $$
\CV\left(\limsup\limits_{k\to\infty} S_{n_k}/{n_k}\ge \overline
\mu \right)=1.$$ This together with (I) suffices to  yield the desired
result (II).

Indeed,  choose $n_k=k^k$ for $k\ge 1.$ Set $\overline
S_n:=\sum\limits_{i=1}^n(X_i-\overline \mu)$, then
$$\begin{array}{lcl}
\displaystyle\CV\left(\frac{ S_{n_k}-S_{n_{k-1}} } {n_k-n_{k-1}}\ge \overline
\mu-\epsilon\right)
&=& \displaystyle\CV\left(\frac{ S_{n_k-n_{k-1}} }
{n_k-n_{k-1}}\ge \overline \mu-\epsilon\right)\\
&=&\displaystyle\CV\left(\frac{ S_{n_k-n_{k-1}}-(n_k-n_{k-1})\overline \mu }
{n_k-n_{k-1}}\ge -\epsilon\right) \\
&=& \displaystyle\CV\left(\frac{ \overline S_{n_k-n_{k-1}}} {n_k-n_{k-1}}\ge
-\epsilon\right)\\
&\ge& \displaystyle\CE[\phi(\frac{ \overline S_{n_k-n_{k-1}}} {n_k-n_{k-1}})],
\end{array}
$$
where $\phi(x)$ is defined by
$$
\phi(x)=\left\{
\begin{array}{l}
1-e^{-(x+\epsilon)},\quad x\ge -\epsilon;\\
0,      \quad\quad\quad\quad  \quad  x<-\epsilon.
\end{array}
\right.
$$
Considering the sequence of IID random variables $\{X_i-\overline
\mu\}_{i=1}^\infty.$
 Obviously $$\CE[X_i-\overline \mu]=0,\quad \Ce[X_i-\overline
\mu]=-(\overline \mu-\underline \mu).$$

Applying Lemma \ref{le4}, we have, $n_k-n_{k-1}\to\infty$ as $k\to \infty$
and
$$\liminf\limits_{n\to\infty}\CE\left[\phi\left(\frac{ \overline S_{n_k-n_{k-1}}} {n_k-n_{k-1}}\right)\right]\ge
\sup_{-(\overline \mu-\underline \mu)\leq y\leq
0}\phi(y)=\phi(0)=1-e^{-\epsilon}>0.$$ Thus
 $$\sum_{k=1}^\infty \CV\left(\frac{ S_{n_k}-S_{n_{k-1}} } {n_k-n_{k-1}}\ge
\overline \mu-\epsilon\right)\ge \sum_{k=1}^\infty\CE\left[\phi\left(\frac{
\overline S_{n_k-n_{k-1}}} {n_k-n_{k-1}}\right)\right]=\infty.$$
From the fact that $\{
 S_{n_k}-S_{n_{k-1}}\}_{k\ge 1}$ is a sequence of independent random variables for $k\ge 1$.
 Using the second Borel-Cantelli Lemma, we have
 $$
\limsup\limits_{k\to\infty}\frac{ S_{n_k}-S_{n_{k-1}} }
{n_k-n_{k-1}}\ge \overline \mu-\epsilon,\quad \hbox{a.s.} \; \CV.
$$
But
$$
\frac{S_{n_k}}{n_k}\ge \frac{ S_{n_k}-S_{n_{k-1}} }
{n_k-n_{k-1}}\cdot\frac{n_k-n_{k-1}}{n_k}-\frac{|S_{n_{k-1}}|}{n_{k-1}}\cdot\frac{{n_{k-1}}}{n_k}.$$
Note the fact that
$$
\frac{n_k-n_{k-1}}{n_k}\to 1, \quad \frac{{n_{k-1}}}{n_k}\to 0, \;
\hbox{as} \; k\to \infty.
$$
and
$$
\limsup\limits_{n\to\infty} S_n/n\leq \overline \mu, \quad
\limsup\limits_{n\to\infty} (-S_n)/n\leq -\underline \mu, \quad \hbox{a.s.}\ v.
$$
we have $$\limsup\limits_{n\to\infty} |S_n|/n\leq\max\{ |\overline
\mu|,|\underline \mu|\},\quad \hbox{a.s.} \ v.$$ Hence,

$$
\limsup\limits_{k\to\infty}\frac{S_{n_k}}{n_k}\ge
\limsup\limits_{k\to\infty}\frac{ S_{n_k}-S_{n_{k-1}}}
{n_k-n_{k-1}}\lim\limits_{k\to\infty}\frac{n_k-n_{k-1}}{n_k}-
\limsup\limits_{k\to\infty}\frac{|S_{n_{k-1}}|}{n_{k-1}}\lim\limits_{k\to\infty}\frac{{n_{k-1}}}{n_k}.
$$
We conclude that
$$
\limsup\limits_{k\to\infty}\frac{S_{n_k}}{n_k}\ge \overline
\mu-\epsilon,\quad \hbox{ a.s. }\; \CV. $$ Since $\epsilon$ is
arbitrary and $\mathbb{V}$ is upper continuous, we have
$$\CV\left(\limsup\limits_{k\to\infty} S_{n_k}/{n_k}\ge
\overline \mu \right)=1.$$
By (I), we know $\CV\left(\limsup\limits_{n\to\infty} S_{n}/{n}>
\overline \mu \right)=0$, thus
$$\begin{array}{lcl}
 \CV\left(\limsup\limits_{n\to\infty} S_{n}/{n}= \overline \mu \right)
&=&\CV\left(\limsup\limits_{n\to\infty} S_{n}/{n}= \overline \mu
\right)+\CV\left(\limsup\limits_{n\to\infty} S_{n}/{n}> \overline
\mu \right)\\
&\ge & \CV\left(\limsup\limits_{n\to\infty} S_{n}/{n}\ge \overline
\mu \right)=1.
\end{array}$$

Considering the sequence of $\{-X_n\}_{n=1}^\infty,$  we have
$$\CV\left(\limsup_{n\to\infty} (-S_{n})/{n}=\CE[-X_1]\right)=1.$$
Therefore,
$$\CV\left(\liminf_{n\to\infty} S_{n}/{n}=-\CE[-X_1]\right)=1.$$
But $\underline \mu=-\CE[-X_1],$ thus
$$\CV\left(\liminf_{n\to\infty} S_{n}/{n}=\underline \mu\right)=1.$$

The proof of (II) is complete.

To prove (III), we only need to prove that, for $b\in (\underline \mu, \overline \mu)$,
$$
\CV\left(\liminf_{n\to\infty}|S_{n}/{n}-b|=0\right)=1.
$$
To do so, we only need to prove that for any $\epsilon>0$ there exists an increasing
subsequence $\{n_k\}$ of $\mathbb{N}$ such that for any
$b\in (\underline \mu, \overline \mu)$,
\begin{equation}\label{R3}
\CV\left(\bigcap_{m=1}^{\infty}\bigcup_{k=m}^{\infty}\{|S_{n_k}/{n_k}-b|\leq\epsilon\}\right)=1.
\end{equation}
Indeed, for any $0<\epsilon \leq \min\{ \overline\mu -b,
b-\underline \mu\},$ let us choose $n_k=k^k$ for $k\ge 1.$

Set $\overline S_n:=\sum\limits_{i=1}^n(X_i-b),$ then
$$\begin{array}{lcl}
\displaystyle\CV\left(\Big|\frac{ S_{n_k}-S_{n_{k-1}} } {n_k-n_{k-1}}-b\Big|\leq
\epsilon\right) &=&\displaystyle \CV\left(\Big|\frac{ S_{n_k-n_{k-1}} }
{n_k-n_{k-1}}-b\Big|\leq \epsilon\right)\\
&=&\displaystyle\CV\left(\Big|\frac{ S_{n_k-n_{k-1}}-(n_k-n_{k-1}) b }
{n_k-n_{k-1}}\Big|\leq \epsilon\right) \\
&=& \displaystyle\CV\left(\Big|\frac{ \overline S_{n_k-n_{k-1}}} {n_k-n_{k-1}}\Big|\leq
\epsilon\right)\\
&\ge&\displaystyle \CE\left[\phi\left(\frac{ \overline S_{n_k-n_{k-1}}} {n_k-n_{k-1}}\right)\right]
\end{array}
$$
where $\phi(x)$ is defined by
$$
\phi(x)=\left\{
\begin{array}{l}
1-e^{|x|-\epsilon},\quad\quad |x|\leq \epsilon;\\
0,      \quad\quad\quad\quad  \quad\;\;  |x|>\epsilon.
\end{array}
\right.
$$
Considering the sequence of IID random variables
$\{X_i-b\}_{i=1}^\infty.$
 Obviously $$\CE[X_i-b]=\overline\mu-b>0,\quad \Ce[X_i-
b]=\underline \mu-b<0.$$
Applying Lemma \ref{le4}, we have
$$\lim\inf\limits_{k\to\infty}\CE\left[\phi\left(\frac{ \overline S_{n_k-n_{k-1}}} {n_k-n_{k-1}}\right)\right]\ge
\sup_{\underline \mu-b\leq y\leq \overline\mu-b
}\phi(y)=\phi(0)=1-e^{-\epsilon}>0.$$ Thus
 $$\sum_{k=1}^\infty \CV\left(\Big|\frac{ S_{n_k}-S_{n_{k-1}} }
 {n_k-n_{k-1}}-b\Big|\leq
\epsilon\right)\geq \sum_{k=1}^\infty\CE\left[\phi\left(\frac{ \overline
S_{n_k-n_{k-1}}} {n_k-n_{k-1}}\right)\right]=\infty.$$

Note the fact that the sequence of
$\{S_{n_k}-S_{n_{k-1}}\}_{k\geq1}$ is independent for all $k\ge 1$ .
Using the second Borel-Cantelli Lemma, we have
 $$
\liminf_{k\to\infty}\Big|\frac{ S_{n_k}-S_{n_{k-1}} }
{n_k-n_{k-1}}-b\Big|\leq \epsilon,\quad \hbox{a.s.} \; \CV.
$$
But
\begin{equation}\label{ED}
\left|\frac{S_{n_k}}{n_k}-b\right|\leq \left|\frac{
S_{n_k}-S_{n_{k-1}}
}{n_k-n_{k-1}}-b\right|\cdot\frac{n_k-n_{k-1}}{n_k}+\left[\frac{|S_{n_{k-1}}|}{n_{k-1}}+|b|\right]\frac{{n_{k-1}}}{n_k}.
\end{equation}
Noting that,
$$
\frac{n_k-n_{k-1}}{n_k}\to 1, \quad \frac{{n_{k-1}}}{n_k}\to 0,\;
\hbox{as} \; k\to \infty
$$
and
$$
\limsup_{n\to\infty} S_n/n\leq \overline \mu, \quad
\limsup_{n\to\infty} (-S_n)/n\leq -\underline \mu,\quad \hbox{a.s.} \; v
$$
which implies
$$
\limsup_{n\to\infty} |S_n|/n\leq \max\{|\overline \mu|, |\underline
\mu|\}<\infty\quad \hbox{a.s.} \; v.
$$
Hence, from inequality(\ref{ED}), for any $\epsilon >0,$
$$
\liminf_{k\to\infty}\Big|\frac{S_{n_k}}{n_k}-b\Big|\leq \epsilon,\quad \hbox{
a.s.}\; \CV.
$$i.e.,
$$\CV\left(\liminf_{n\to \infty} |S_{n}/{n}-b|\leq \epsilon\right)=1.$$
Since $\epsilon$ is arbitrary, we have

$$\CV\left(\liminf_{n\to \infty} |S_{n}/{n}-b|=0\right)=1.$$

 The proof of (III) is complete.
\end{proof}

\Acknowledgements{This work was supported by National Natural Science Foundation of China (Grant No. 11231005). The author thanks the anonymous referees for their valuable comments which greatly improved this paper.}

\end{document}